\documentclass[12pt]{article}
\usepackage{amssymb,t1enc}
\newcommand{\Q}{\mathbb Q}
\newcommand{\R}{\mathbb R}
\newcommand{\C}{\mathbb C}

\begin{document}
\thispagestyle{empty}
\par
\noindent
\centerline{{\large Some conjectures on addition and multiplication of complex (real) numbers}}
\vskip 0.2truecm
\par
\centerline{{\large Apoloniusz Tyszka}}
\vskip 0.2truecm
\par
\noindent
{\bf Abstract.} We discuss conjectures related to the following two conjectures:
\par
\noindent
{\bf (I)} (see \cite{Tyszka2}) for each complex numbers $x_1,\ldots,x_n$
there exist rationals $y_1,\ldots,y_n \in [-2^{n-1},2^{n-1}]$ such that
\begin{equation}
\forall i \in \{1,\ldots,n\} ~(x_i=1 \Rightarrow y_i=1)
\end{equation}
\begin{equation}
\forall i,j,k \in \{1,\ldots,n\} ~(x_i+x_j=x_k \Rightarrow y_i+y_j=y_k)
\end{equation}
\par
\noindent
{\bf (II)} (see \cite{Tyszka1},~\cite{Tyszka2}) for each complex (real) numbers
$x_1,\ldots,x_n$ there exist complex (real) numbers $y_1,\ldots,y_n$ such that
\begin{equation}
\forall i \in \{1,\ldots,n\} ~|y_i| \leq 2^{\textstyle 2^{n-2}}
\end{equation}
\begin{equation}
\forall i \in \{1,\ldots,n\} ~(x_i=1 \Rightarrow y_i=1)
\end{equation}
\begin{equation}
\forall i,j,k \in \{1,\ldots,n\} ~(x_i+x_j=x_k \Rightarrow y_i+y_j=y_k)
\end{equation}
\begin{equation}
\forall i,j,k \in \{1,\ldots,n\} ~(x_i \cdot x_j=x_k \Rightarrow y_i \cdot y_j=y_k)
\end{equation}
\vskip 1.0truecm
\par
For a positive integer $n$ we define the set of equations $W_n$ by
\par
\noindent
\centerline{$W_n=\{x_i=1:~1 \leq i \leq n\} \cup \{x_i+x_j=x_k:~1 \leq i \leq j \leq n,~ 1 \leq k \leq n\}$}
\vskip 0.2truecm
\par
\noindent
Let $S \subseteq W_n$ be a system consistent over $\C$.
Then $S$ has a solution which consists of rationals belonging to
$[-(\sqrt{5})^{n-1},(\sqrt{5})^{n-1}]$, see [7,~Theorem~9].
Conjecture {\bf (I)} states that $S$ has a solution which consists of rationals
belonging to $[-2^{n-1},2^{n-1}]$.
\vskip 0.2truecm
\par
Concerning Conjecture {\bf (I)}, estimation by $2^{n-1}$
is the best estimation. Indeed, the system
\begin{displaymath}
\left\{
\begin{array}{rcl}
x_1&=& 1 \\
x_1 + x_1 &=& x_2 \\
x_2 + x_2 &=& x_3 \\
x_3 + x_3 &=& x_4 \\
&...& \\
x_{n-1} + x_{n-1} &=&x_n
\end{array}
\right.
\end{displaymath}
has precisely one complex solution: $(1,2,4,8,\ldots,2^{n-2},2^{n-1})$.
Concerning \hbox{Conjecture {\bf (II)},}
for $n=1$ estimation by $2^{\textstyle 2^{n-2}}$ can be replaced
by estimation by $1$. For $n>1$ estimation by $2^{\textstyle 2^{n-2}}$
is the best estimation, see \cite{Tyszka2}.
\vskip 0.2truecm
\par
For each consistent system $S \subseteq W_n$ there exists $J \subseteq \{1,\ldots,n\}$
such that the system $S \cup \{x_i+x_i=x_i:~i \in J\}$ has a unique solution $(x_1,\ldots,x_n)$,
see the proof of \hbox{Theorem 9} in \cite{Tyszka2}.
If any $S \subseteq W_n$ has a unique solution $(x_1,\ldots,x_n)$,
then by Cramer's rule each $x_i$ is a quotient of two determinants.
Since these determinants have entries among $-1$, $0$, $1$, $2$, each $x_i$ is rational.
\vskip 0.2truecm
\par
For proving Conjecture {\bf (I)}, without loss of generality we can assume that the equation $x_1=1$ belongs to $S$
and all equations $x_i=1$ ($i>1$) do not belong to $S$.
Indeed, if $i>1$ and the equation $x_i=1$ belongs to $S$,
then we replace $x_i$ by $x_1$ in all equations belonging to $S$.
In this way the problem reduces to the same problem
with a smaller number of variables, for details see the text after Conjecture 3.
Therefore, for proving \hbox{Conjecture {\bf (I)}}
it is sufficient to consider only these systems $S \subseteq W_n$ of $n$ equations which satisfy the following conditions:
\par
$S$ contains the equation $x_1=1$ and $n-1$ equations of the form $x_i+x_j=x_k$ ($i,j,k \in \{1,\ldots,n\}$),
\par
joining the equation $x_1=0$ and the aforementioned $n-1$ equations of the form $x_i+x_j=x_k$,
we get a system of $n$ linearly independent equations.
\vskip 0.2truecm
\par
By the observations from the last two paragraphs, the following code in {\sl MuPAD}
yields a probabilistic confirmation of Conjecture {\bf (I)}.
The value of $n$ is set, for example, to $5$.
The number of iterations is set, for example, to $1000$.
We use the algorithm which terminates with probability $1$.
For another algorithm, implemented in {\sl Mathematica}, see \cite{Kozlowski}.
\begin{quote}
\begin{verbatim}
SEED:=time():
r:=random(1..5):
idmatrix:=matrix::identity(5):
u:=linalg::row(idmatrix,i) $i=1..5:
max_norm:=1:
for k from 1 to 1000 do
a:=linalg::row(idmatrix,1):
rank:=1:
while rank<5 do
m:=matrix(u[r()])+matrix(u[r()])-matrix(u[r()]):
a1:=linalg::stackMatrix(a,m):
rank1:=linalg::rank(a1):
if rank1 > rank then a:=linalg::stackMatrix(a,m) end_if:
rank:=linalg::rank(a):
end_while:
x:=(a^-1)*linalg::col(idmatrix,1):
max_norm:=max(max_norm,norm(x)):
print(max_norm):
end_for:
\end{verbatim}
\end{quote}
\par
Conjecture {\bf (I)} holds true for each $n \leq 4$.
It follows from the following \hbox{Observation 1.}
\vskip 0.2truecm
\par
\noindent
{\bf Observation 1} ([7,~p.~23]). If $n \leq 4$ and $(x_1,\ldots,x_n) \in {\C}^n$ solves $S$,
then \hbox{some $(\widehat{x_1},\ldots,\widehat{x_n})$} solves $S$, where each $\widehat{x_i}$ is suitably
chosen from $\{x_i,0,1,2,\frac{1}{2}\} \cap \{r \in \Q:~|r| \leq 2^{n-1}\}$.
\vskip 0.2truecm
\par
Let ${\bf A}{\bf x}={\bf b}$ be the matrix representation of the system $S$,
and let ${\bf A}^{\dagger}$ denote Moore-Penrose pseudoinverse of ${\bf A}$.
The system $S$ has a unique solution ${\bf x_0}$ with minimal Euclidean norm, and this
element is given by ${\bf x_0}={\bf A}^{\dagger}{\bf b}$, see [5,~p.~423].
\vskip 0.2truecm
\par
For any system $S \subseteq W_n$, a vector ${\bf x} \in {\C}^n$ is said to be
a least-squares solution if ${\bf x}$ minimizes the Euclidean norm of ${\bf Ax}-{\bf b}$, see [1,~p.~104].
It is known that ${\bf x_0}={\bf A}^{\dagger}{\bf b}$ is a unique least-squares solution with minimal
Euclidean norm, see [1,~p.~109].
\vskip 0.2truecm
\par
Since ${\bf A}$ has rational entries (the entries are among $-1$, $0$, $1$, $2$),
${\bf A}^{\dagger}$ has also rational entries, see [2,~p.~69] and [3,~p.~193].
Since ${\bf b}$ has rational entries (the entries are among $0$ and~$1$),
${\bf x_0}={\bf A}^{\dagger}{\bf b}$ consists of rationals.
\vskip 0.2truecm
\par
\noindent
{\bf Conjecture 1.} The solution (The least-squares solution) ${\bf x_0}$ consists of numbers belonging
to $[-2^{n-1},2^{n-1}]$.
\vskip 0.2truecm
\par
Conjecture 1 restricted to the case when ${\rm card~} S \leq n$ implies Conjecture {\bf (I)}.
The following code in {\sl MuPAD} yields a probabilistic confirmation of Conjecture 1 restricted
to the case when ${\rm card~} S \leq n$. The value of $n$ is set, for example, to $5$.
The number of iterations is set, for example, to $1000$.
\begin{quote}
\begin{verbatim}
SEED:=time():
r:=random(1..5):
idmatrix:=matrix::identity(5):
u:=linalg::row(idmatrix,i) $i=1..5:
max_norm:=1:
for k from 1 to 1000 do
b:=[1]:
c:=linalg::row(idmatrix,1):
for w from 1 to 4 do
h:=0:
h1:=r():
h2:=r():
h3:=r():
m1:=matrix(u[h1]):
m2:=matrix(u[h2]):
m3:=matrix(u[h3]):
m:=m1+m2-m3:
c:=linalg::stackMatrix(c,m):
if h3=h2 then h:=1 end_if:
b:=append(b,h):
a:=linalg::pseudoInverse(c):
x:=a*matrix(b):
max_norm:=max(max_norm,norm(x)):
end_for:
print(max_norm):
end_for:
\end{verbatim}
\end{quote}
\par
The following Conjecture 2 implies Conjecture {\bf (I)}, see \cite{Tyszka2}.
\vskip 0.2truecm
\par
\noindent
{\bf Conjecture 2} (\cite{Tyszka2}). Let ${\bf B}$ be a matrix with $n-1$ rows and $n$ columns, $n \geq 2$.
Assume that each row of ${\bf B}$, after deleting all zeros, forms a sequence belonging to
\par
\noindent
\centerline{$\{
\langle 1 \rangle,
\langle -1,2 \rangle,
\langle 2,-1 \rangle,
\langle -1,1,1 \rangle,
\langle 1,-1,1 \rangle,
\langle 1,1,-1 \rangle\}$}
\par
\noindent
We conjecture that after deleting any column of ${\bf B}$ we get the matrix
whose determinant has absolute value less than or equal to $2^{n-1}$.
\vskip 0.2truecm
\par
\noindent
{\bf Conjecture 3.} If a system $S \subseteq W_n$ has a unique solution $(x_1,\ldots,x_n)$,
then this solution consists of rationals whose nominators and denominators belong to $[-2^{n-1},2^{n-1}]$.
\vskip 0.2truecm
\par
Conjecture 3 implies Conjecture {\bf (I)}.
The {\sl MuPAD} code below confirms \hbox{Conjecture 3} probabilistically.
As previously, the value of $n$ is set to $5$, the number of iterations is set to $1000$.
We declare that
\vskip 0.2truecm
\par
\noindent
\centerline{$\{i \in \{1,2,3,4,5\}:{\rm ~the~equation~}x_i=1 {\rm ~belongs~to~} S\}=\{1\}$,}
\vskip 0.2truecm
\par
\noindent
but this does not decrease the generality. Indeed, $(0,\ldots,0)$ solves $S$ if all equations
$x_i=1$ do not belong to $S$. In other cases, let
\vskip 0.2truecm
\par
\noindent
\centerline{$I=\{i \in \{1,\ldots,n\}:{\rm~the~equation~} x_i=1 {\rm ~belongs~to~} S$\},}
\vskip 0.2truecm
\par
\noindent
and let $i={\rm min}(I)$. For each $j \in I$ we replace $x_j$ by $x_i$ in all equations belonging
\hbox{to $S$.} We obtain an equivalent system $\widehat{S}$ with $n-{\rm card}(I)+1$ variables.
The system $\widehat{S}$ has a unique solution $(t_1,\ldots,t_{n-{\rm card}(I)+1})$,
and the equation $x_j=1$ belongs to $\widehat{S}$ if and only if $j=i$.
By permuting variables, we may assume that $i=1$.
\begin{quote}
\begin{verbatim}
SEED:=time():
r:=random(1..5):
idmatrix:=matrix::identity(5):
u:=linalg::row(idmatrix,i) $i=1..5:
abs_numer_denom:=[1]:
for k from 1 to 1000 do
c:=linalg::row(idmatrix,1):
rank:=1:
while rank<5 do
m:=matrix(u[r()])+matrix(u[r()])-matrix(u[r()]):
if linalg::rank(linalg::stackMatrix(c,m))>rank
then c:=linalg::stackMatrix(c,m) end_if:
rank:=linalg::rank(c):
end_while:
a:=(c^-1)*linalg::col(idmatrix,1):
for n from 2 to 5 do
abs_numer_denom:=append(abs_numer_denom,abs(numer(a[n]))):
abs_numer_denom:=append(abs_numer_denom,abs(denom(a[n]))):
end_for:
abs_numer_denom:=listlib::removeDuplicates(abs_numer_denom):
print(max(abs_numer_denom)):
end_for:
\end{verbatim}
\end{quote}
\par
The {\sl MuPAD} code below completely confirms Conjecture 3 for $n=5$. We declare that
\par
\noindent
\centerline{$\{i \in \{1,2,3,4,5\}:{\rm ~the~equation~}x_i=1 {\rm ~belongs~to~} S\}=\{1\}$,}
\par
\noindent
but this does not decrease the generality.
\begin{quote}
\begin{verbatim}
p:=[]:
idmatrix:=matrix::identity(5):
u:=linalg::row(idmatrix,i) $i=1..5:
for r1 from 1 to 5 do
for r2 from 1 to 5 do
for r3 from 1 to 5 do
m1:=matrix(u[r1]):
m2:=matrix(u[r2]):
m3:=matrix(u[r3]):
m:=m1+m2-m3:
p:=append(p,m):
end_for:
end_for:
end_for:
p:=listlib::removeDuplicates(p):
p:=listlib::setDifference(p,[linalg::row(idmatrix,1)]):
abs_numer_denom:=[]:
s1:=nops(p)-1:
s2:=nops(p)-2:
s3:=nops(p)-3:
for n3 from 1 to s3 do
w2:=n3+1:
for n2 from w2 to s2 do
w1:=n2+1:
for n1 from w1 to s1 do
w0:=n1+1:
for n0 from w0 to nops(p) do
c3:=linalg::stackMatrix(linalg::row(idmatrix,1),p[n3]):
c2:=linalg::stackMatrix(c3,p[n2]):
c1:=linalg::stackMatrix(c2,p[n1]):
c:=linalg::stackMatrix(c1,p[n0]):
if linalg::rank(c)=5 then
a:=(c^-1)*linalg::col(idmatrix,1):
for n from 2 to 5 do
abs_numer_denom:=append(abs_numer_denom,abs(numer(a[n]))):
abs_numer_denom:=append(abs_numer_denom,abs(denom(a[n]))):
end_for:
abs_numer_denom:=listlib::removeDuplicates(abs_numer_denom):
end_if:
end_for:
end_for:
abs_numer_denom:=sort(abs_numer_denom):
print(abs_numer_denom):
end_for:
end_for:
\end{verbatim}
\end{quote}
\par
The following Conjecture 4 implies Conjecture {\bf (I)}, because each consistent system $S \subseteq W_n$
can be enlarged to a system $\widetilde{S} \subseteq W_n$ with a unique solution $(x_1,\ldots,x_n)$ and
$(x_1,\ldots,x_n) \in {\Q}^n$.
\vskip 0.2truecm
\par
\noindent
{\bf Conjecture 4.} Let rationals $x_1,\ldots,x_n$ satisfy $|x_1| \leq |x_2| \leq \ldots \leq |x_n|$,
and for each $y_1,\ldots,y_n \in \Q$, if
\par
\noindent
\centerline{$\forall i \in \{1,\ldots,n\} ~(x_i=1 \Rightarrow y_i=1)$}
\par
\noindent
and
\par
\noindent
\centerline{$\forall i,j,k \in \{1,\ldots,n\} ~(x_i+x_j=x_k \Rightarrow y_i+y_j=y_k)$,}
\par
\noindent
then $(x_1,\ldots,x_n)=(y_1,\ldots,y_n)$.
We conjecture that for each $i \in \{1,\ldots,n-1\}$ the inequality $|x_i| \geq 1$ implies
$|x_{i+1}| \leq 2 \cdot |x_i|$.
\vskip 0.2truecm
\par
Concerning Conjecture 4, the assumption $|x_i| \geq 1$ is necessary to state that
$|x_{i+1}| \leq 2 \cdot |x_i|$. As a trivial counterexample we have $(0,1)$,
$(\frac{1}{4}, \frac{3}{4}, 1, \frac{3}{2}, 2, 3)$
is a counterexample which consists of positive rationals alone.
The {\sl MuPAD} code below probabilistically confirms Conjecture 4.
The value of $n$ is set, for example, to $5$.
The number of iterations is set, for example, to $1000$.
\begin{quote}
\begin{verbatim}
SEED:=time():
r:=random(1..5):
idmatrix:=matrix::identity(5):
u:=linalg::row(idmatrix,i) $i=1..5:
max_ratio:=1:
for k from 1 to 1000 do
a:=linalg::row(idmatrix,1):
rank:=1:
while rank<5 do
m:=matrix(u[r()])+matrix(u[r()])-matrix(u[r()]):
a1:=linalg::stackMatrix(a,m):
rank1:=linalg::rank(a1):
if rank1 > rank then a:=linalg::stackMatrix(a,m) end_if:
rank:=linalg::rank(a):
end_while:
x:=(a^-1)*linalg::col(idmatrix,1):
xx:=[max(1,abs(x[i])) $i=1..5]:
xxx:=sort(xx):
maxratio:=max([xxx[i+1]/xxx[i] $i=1..4]):
max_ratio:=max(max_ratio,maxratio):
print(max_ratio):
end_for:
\end{verbatim}
\end{quote}
\par
For a positive integer $n$ we define the set of equations $E_n$ by
\par
\noindent
\centerline{$E_n=\{x_i=1:~1 \leq i \leq n\}~\cup$}
\par
\noindent
\centerline{$\{x_i+x_j=x_k:~1 \leq i \leq j \leq n,~1 \leq k \leq n\}
\cup
\{x_i \cdot x_j=x_k:~1 \leq i \leq j \leq n,~1 \leq k \leq n\}$}
\vskip 0.2truecm
\par
\noindent
Let $T \subseteq E_n$ be a system consistent over $\C$ (over $\R$).
Conjecture {\bf (II)} states that $T$ has a complex (real) solution
which consists of numbers whose absolute values belong to $[0,2^{\textstyle 2^{n-2}}]$.
Both for complex and real case, we conjecture that each solution of $T$ with minimal
Euclidean norm consists of numbers whose absolute values belong to $[0,2^{\textstyle 2^{n-2}}]$.
This conjecture implies Conjecture {\bf (II)}. Conjecture {\bf (II)} holds true for each $n \leq 4$.
It follows from the following Observation 2.
\vskip 0.2truecm
\par
\noindent
{\bf Observation 2}~([7,~p.~7]). If $n \leq 4$ and $(x_1,\ldots,x_n) \in {\C}^n~({\R}^n)$ solves $T$,
then some $(\widehat{x_1},\ldots,\widehat{x_n})$ solves $T$, where each $\widehat{x_i}$ is suitably chosen from
$\{x_i,0,1,2,\frac{1}{2}\} \cap$ \hbox{$\{z \in \C~(\R):~|z| \leq 2^{\textstyle 2^{n-2}}\}$.}
\vskip 0.2truecm
\par
Let us consider the following four conjectures, analogical conjectures seem to be true for $\R$.
\begin{description}
\item{{\bf (5a)}}
If a system $S \subseteq E_n$ is consistent over $\C$ and maximal with respect to inclusion,
then each solution of $S$ belongs to\\
$\{(x_1,\ldots,x_n) \in {\C}^n:~|x_1| \leq 2^{\textstyle 2^{n-2}}~\wedge~\ldots~\wedge~|x_n| \leq 2^{\textstyle 2^{n-2}}\}$.
\item{{\bf (5b)}}
If a system $S \subseteq E_n$ is consistent over $\C$ and maximal with respect to inclusion,
then $S$ has a finite number of solutions $(x_1,\ldots,x_n)$.
\item{{\bf (5c)}}
If the equation $x_1=1$ belongs to $S \subseteq E_n$ and $S$ has a finite number of complex solutions $(x_1,\ldots,x_n)$,
then each such solution belongs to\\
$\{(x_1,\ldots,x_n) \in {\C}^n:~|x_1| \leq 2^{\textstyle 2^{n-2}}~\wedge~\ldots~\wedge~|x_n| \leq 2^{\textstyle 2^{n-2}}\}$.
\item{{\bf (5d)}}
If a system $S \subseteq E_n$ has a finite number of complex solutions $(x_1,\ldots,x_n)$,
then each such solution belongs to\\
$\{(x_1,\ldots,x_n) \in {\C}^n:~|x_1| \leq 2^{\textstyle 2^{n-1}}~\wedge~\ldots~\wedge~|x_n| \leq 2^{\textstyle 2^{n-1}}\}$.
\end{description}
\par
Conjecture 5a strengthens Conjecture {\bf (II)} for $\C$. The conjunction of Conjectures 5b and 5c implies Conjecture 5a.
\vskip 0.2truecm
\par
Concerning Conjecture 5d, for $n=1$ estimation by $2^{\textstyle 2^{n-1}}$ can be replaced
by estimation by $1$. For $n>1$ estimation by $2^{\textstyle 2^{n-1}}$ is the best estimation. Indeed, the system
\begin{displaymath}
\left\{
\begin{array}{rcl}
x_1+x_1 &=& x_2 \\
x_1 \cdot x_1 &=& x_2 \\
x_2 \cdot x_2 &=& x_3 \\
x_3 \cdot x_3 &=& x_4 \\
&...& \\
x_{n-1} \cdot x_{n-1} &=&x_n
\end{array}
\right.
\end{displaymath}
has precisely two complex solutions, $(0,\ldots,0)$, $(2,4,16,256,\ldots,2^{\textstyle 2^{n-2}},2^{\textstyle 2^{n-1}})$.
\vskip 0.2truecm
\par
The following code in {\sl MuPAD} yields a probabilistic confirmation of Conjectures 5b and 5c.
The value of $n$ is set, for example, to $5$. The number of iterations is set, for example, to $1000$.
\begin{quote}
\begin{verbatim}
SEED:=time():
p:=[v-1,x-1,y-1,z-1]:
var:=[1,v,x,y,z]:
for i from 1 to 5 do
for j from i to 5 do
for k from 1 to 5 do
p:=append(p,var[i]+var[j]-var[k]):
p:=append(p,var[i]*var[j]-var[k]):
end_for:
end_for:
end_for:
p:=listlib::removeDuplicates(p):
max_abs_value:=1:
for r from 1 to 1000 do
q:=combinat::permutations::random(p):
syst:=[t-v-x-y-z]:
w:=1:
repeat
if groebner::dimension(append(syst,q[w]))>-1
then syst:=append(syst,q[w]) end_if:
w:=w+1:
until (groebner::dimension(syst)=0 or w>nops(q)) end:
d:=groebner::dimension(syst):
if d>0 then print("Conjecture 5b is false") end_if:
if d=0 then
sol:=numeric::solve(syst):
for m from 1 to nops(sol) do
for n from 2 to 5 do
max_abs_value:=max(max_abs_value,abs(sol[m][n][2])):
end_for:
end_for:
end_if:
print(max_abs_value);
end_for:
\end{verbatim}
\end{quote}
\par
\noindent
If we replace
\begin{quote}
\begin{verbatim}
p:=[v-1,x-1,y-1,z-1]:    by    p:=[]:
var:=[1,v,x,y,z]:        by    var:=[u,v,x,y,z]:
max_abs_value:=1:        by    max_abs_value:=0:
syst:=[t-v-x-y-z]:       by    syst:=[t-u-v-x-y-z]:
for n from 2 to 5 do     by    for n from 2 to 6 do
\end{verbatim}
\end{quote}
then we get a code for a probabilistic confirmation of Conjecture 5d.
\vskip 0.2truecm
\par
It seems that for each integers $x_1,\ldots,x_n$ there exist integers
$y_1,\ldots,y_n \in [-2^{n-1},2^{n-1}]$ with properties (1) and (2), cf.~[7,~Theorem 10].
However, not for each integers $x_1,\ldots,x_n$ there exist integers $y_1,\ldots,y_n$
with properties (3)-(6), see [7,~pp.~15--16].
\vskip 0.2truecm
\par
The author used MuPAD Pro 4.0.6.
SciFace Software GmbH \& Co. KG, the maker of MuPAD Pro, has been acquired by The MathWorks,
the maker of MATLAB technical computing software. MuPAD Pro is no longer sold as a standalone product.
A new package, Symbolic Math Toolbox 5.1 requires MATLAB and provides large compatibility with existing
MuPAD Pro applications.

Apoloniusz Tyszka\\
Technical Faculty\\
Hugo Ko\l{}\l{}\k{a}taj University\\
Balicka 116B, 30-149 Krak\'ow, Poland\\
E-mail: {\it rttyszka@cyf-kr.edu.pl}
\end{document}